\documentclass{article}
\usepackage{amssymb,latexsym,amscd}

\newcommand{\se}[1]{{\section{#1}} {\setcounter{equation}{0}}}
\newtheorem{theorem}{Theorem}[section]
\newtheorem{lm}{Lemma}[section]

\newtheorem{de}{Definition}[section]
\newtheorem{co}{Corollary}[section]

\def\k{{K\"{a}hler }}

\def\cy{{Calabi-Yau }}

\begin{document}
\hbadness=10000
\title{{\bf Generalized special Lagrangian torus fibration for Calabi-Yau hypersurfaces in toric varieties I}}
\author{Wei-Dong Ruan\\
Department of Mathematics\\
University of Illinois at Chicago\\
Chicago, IL 60607\\}
\footnotetext{Partially supported by NSF Grant DMS-0104150.}
\maketitle
\begin{abstract}
In this paper we start the program of constructing generalized special Lagrangian torus fibrations for Calabi-Yau hypersurfaces in toric varieties near the large complex limit, with respect to the restriction of a toric metric on the toric variety to the Calabi-Yau hypersurface. The construction is based on the deformation of the standard toric generalized special Lagrangian torus fibration of the large complex limit $X_0$. In this paper, we will deal with the region near the smooth top dimensional torus fibres of $X_0$ and its mirror dual situation: the region near the 0-dimensional fibres of $X_0$.
\end{abstract}
\se{Introduction}
Let $(X,\omega_g,\Omega)$ be a Calabi-Yau manifold, where $\omega_g$ is the \k form of a \k metric $g$ (not necessarily Calabi-Yau) and $\Omega$ is a holomorphic volume form. A middle dimensional real submanifold $L\subset X$ is called a {\bf generalized special Lagrangian submanifold} if $L$ is Lagrangian with respect to the \k form $\omega_g$ and $\Omega$ has constant phase along $L$. Namely, the phase angle function $\theta_1$ is constant, where $\Omega|_L = e^{i\theta} d{\rm Vol}_g$, $\theta = \theta_1 + i\theta_2$. Notice that $\theta_2$ is the restriction on $L$ of a global function $\theta_2$ on $X$ satisfying $\Omega \wedge \bar{\Omega} = e^{-2\theta_2} \omega_g^n$. When $\theta_2$ is a constant, $g$ is a Calabi-Yau metric (\cite{yau}) and $L$ is called {\bf special Lagrangian submanifold} of $X$. (More background and details on these concepts can be found in \cite{Joyce1}.) In this paper, we will mainly consider generalized special Lagrangian submanifolds and generalized special Lagrangian torus fibration for \cy hypersurfaces in toric varieties near large complex limit.\\

The concept of special Lagrangian submanifold was first introduced in the seminal work of Harvey and Lawson \cite{HL} in the more general context of calibrated geometry. In \cite{HL} some important examples were constructed. After \cite{HL}, the most notable progress on special Lagrangian was made in the work of Mclean \cite{ML}. Mclean's work was actually on the more general framework of deformation of calibrated submanifolds. Specifying to the case of special Lagrangian, he proved that the local deformation of a special Lagrangian submanifold $L$ is unobstructed and the local moduli has real dimension $h^1(L)$. This important work did not attract much attention and was not published until the proposal of SYZ conjecture \cite{SYZ} in mirror symmetry. \\

According to SYZ proposal, on each Calabi-Yau manifold there should be a special Lagrangian torus fibration. This conjectural special Lagrangian torus fibration structure of Calabi-Yau manifolds is further used to give a possible geometric explanation of the mirror symmetry conjecture. SYZ conjecture brought a lot of new attention to the study of special Lagrangian, and more generally calibrated geometry. \\

In complex dimension 2, the SYZ special Lagrangian fibrations can be reduced to the classical elliptic fibrations for K3 surfaces under the hyper\k twist, which generically have 24 simple rational nodal singular fibres. The detailed knowledge of the two-dimensional SYZ fibrations naturally would suggest that SYZ fibration in dimension three should be smooth with codimension 2 singular locus. This had indeed been the prevailing view for quite a while. However evidences later indicated that this prevailing view might not be true. One evidence appeared in my construction of Lagrangian torus fibrations for Calabi-Yau hypersurfaces in toric varieties via gradient flow method (\cite{lag1, lag2, lag3, toric}). The singular loci for the fibrations produced by the natural Hamiltonian-gradient flow are of amoeba nature and of codimension 1! Incidentally, the fibration maps are not smooth! This made me suspect that the actual special Lagrangian fibrations might also exhibit similar behavior (see the remark in the last section of \cite{toric}). Even with this evidence, people in favor of the original point of view could still argue that special Lagrangian condition might miraculously turn the fibration maps into smooth maps and make singular locus codimension 2. A more convincing evidence came from the important local example of special Lagrangian fibration with codimension 1 singular locus by D. Joyce \cite{Joyce}. His fibration map is also not smooth. After these evidences, the current prevailing view is that the special Lagrangian torus fibrations for Calabi-Yau manifolds (if exist) should be non-smooth maps and have codimension 1 singular loci, very similar to the Lagrangian fibrations constructed via gradient flow method in \cite{lag1, lag2, lag3, toric}. Nevertheless, up till now there is still not a single global construction of (generalized) special Lagrangian torus fibration with codimension 1 singular locus for a compact Calabi-Yau 3-fold to the author's knowledge. \\

After SYZ conjecture, especially after the examples of D. Joyce, there has been an explosion of work on singular special Lagrangians. We will not discuss this development in detail, since our work at this stage only concerns smooth (generalized) special Lagrangian torus. Interested readers may consult the series of papers on this subject by D. Joyce (starting with \cite{Joyce}) and the references therein. \\

Our current work concentrates on the global structure of (generalized) special Lagrangian fibrations and can be viewed as an attempt toward verifying my conjecture (\cite{lag2}) that the actual special Lagrangian fibrations should resemble the structure of Lagrangian fibrations constructed via gradient flow. Instead of constructing special Lagrangian submanifold by hand out of nowhere, which is very hard indeed, we look closer to the physics origin of the problem and notice that the explicit toric Lagrangian fibration of the toric large complex limit, based on which we produced Lagrangian fibrations for Calabi-Yau hypersurfaces near the large complex limit via gradient flow approach, is actually also a generalized special Lagrangian toric fibration over a sphere! Our method is a combination of the Hamilton-gradient flow we used in \cite{lag1, lag2, lag3, toric} and the deformation of (generalized) special Lagrangian manifolds similar to Mclean's approach using some implicit function theorem.\\ 

More precisely, the fibration we construct is a generalized special Lagrangian fibration of an open set of the Calabi-Yau hypersurface over an open set of the base sphere with smooth torus fibres. The first fact that makes such (partial) fibration interesting is that (generalized) special Lagrangian fibration is rigid and free of local Hamiltonian deformation. One can be sure that any (generalized) special Lagrangian fibre one constructs will remain to be part of the ultimate global (generalized) special Lagrangian fibration for the Calabi-Yau. This suggests that one can construct the (generalized) special Lagrangian fibration piece by piece and these pieces will automatically coincide in the overlap and form a global fibration. In comparison, such (partial) fibration will be of far less interest for Lagrangian fibrations precisely due to the local Hamiltonian freedom of Lagrangian fibrations. Secondly, our fibration will be over an open set of the base sphere whose compliment is a fattening of the singular locus of the corresponding Lagrangian fibration constructed in \cite{lag1, lag2, lag3, toric}. This enables us to verify all the non-trivial monodromies and they coincide with the Lagrangian case, which is a good indication that the ultimate global (generalized) special Lagrangian fibration for the Calabi-Yau hypersurface might resemble the structure of Lagrangian fibration constructed via gradient flow. For this reason, we call our (partial) fibration {\it monodromy representing}. In a certain sense, the construction of monodromy representing fibration is the construction of smooth fibres, which nicely compliments the works of Joyce and others on singular special Lagrangians. Hopefully one can eventually fill in the singular fibres and their neighborhood and achieve the construction of global (generalized) special Lagrangian fibration.\\

{\bf Remark:} Although the SYZ mirror conjecture is formulated in terms of special Lagrangian fibrations, generalized special Lagrangian fibration is in a certain sense a close enough substitute. Firstly, generalized special Lagrangian submanifolds are minimal submanifolds with respect to the conformal \k metric $\tilde{g} = e^{-\frac{2}{n}\theta_2} g$, just like the special Lagrangian submanifolds. So they could exhibit similar singularity structures. Secondly, as for special Lagrangian submanifolds, the deformation space of generalized special Lagrangian submanifold $L$ in $X$ also has dimension $h_1(L)$. On the other hand, generalized special Lagrangian submanifolds have the advantage that they make sense for any \k metric, so we do not need to work with the Calabi-Yau metric (\cite{yau}), whose structure is very mysterious to this day.\\

\se{Overview of the construction}
Let $(P_\Delta,\omega)$ be a toric variety whose moment map image (with respect to the toric \k form $\omega$) is the real convex polyhedron $\Delta \subset M_{\mathbb{R}}$. Also assume that the anti-canonical class of $P_\Delta$ is represented by an integral reflexive convex polyhedron $\Delta_0\subset M$ and the unique interior point of $\Delta_0$ is the origin of the lattice $M$. Integral points $m \in \Delta_0$ correspond to holomorphic toric sections $s_m$ of the anti-canonical bundle. For the unique interior point $m_o$ of $\Delta_0$, $s_{m_o}$ is the section of the anti-canonical bundle that vanishes to first order along each toric divisor of $P_\Delta$.\\

Let $\{w_m\}_{m\in \Delta_0}$ be a strictly convex function on $\Delta_0$ such that $w_m> 0$ for $m\in \Delta_0\setminus \{m_o\}$ and $w_{m_o}\ll 0$. Define

\[
\tilde{s}_t = s_{m_o} + ts,\ \ s = \sum_{m\in \Delta_0\setminus \{m_o\}} a_m s_m,\  {\rm where}\ |a_m| = \tau^{w_m},\ {\rm for}\ m\in \Delta_0\setminus \{m_o\}.
\]

Let $X_t = \{\tilde{s}_t^{-1}(0)\}$. Then $\{X_t\}$ is a 1-parameter family of Calabi-Yau hypersurfaces. $X_0 = \{s_{m_o}^{-1}(0)\}$ is the so-called large complex limit. $X_t$ is said to be near the large complex limit if $\tau$ and $t$ are small and $t \leq \tau^{-w_{m_o}}$.\\

Since $X_0$ is toric, the moment map induces the standard generalized special Lagrangian fibration $\hat{\pi}_0: X_0 \rightarrow \partial \Delta$ with respect to the toric holomorphic volume form. In \cite{lag1, lag2, lag3, toric}, we constructed Lagrangian torus fibration for $X_t$ when $X_t$ is near the large complex limit, using the Hamiltonian-gradient flow to deform this fibration for $X_0$ symplectically to the desired Lagrangian fibration $\hat{\pi}_t: X_t \rightarrow \partial \Delta$ for such $X_t$. The (topological) singular set of the fibration map $\hat{\pi}_t$ is $C= X_t \cap {\rm Sing}(X_0)$, which is independent of $t$. The corresponding singular locus $\tilde{\Gamma} = \hat{\pi}_0(C)$ is also independent of $t$. When $X_t$ is near the large complex limit, $\tilde{\Gamma} \subset \partial \Delta$ exhibits amoeba structure that is a fattening of a graph $\Gamma \subset \partial \Delta$. It was conjectured in \cite{lag2} that the singular locus for the (generalized) special Lagrangian torus fibration should resemble the singular locus $\tilde{\Gamma}$ of the Lagrangian torus fibration.\\ 

In this work we will use similar idea to construct generalized special Lagrangian fibration for $X_t$ by deformation from the standard fibration of the large complex limit $X_0$.\\

At a first glance, the construction of special Lagrangian fibration (middle dimensional minimal surface PDE) appears to be much harder than the construction of Lagrangian fibration (with Hamiltonian deformation freedom). However, since special Lagrangian submanifold is rigid, the generalized special Lagrangian fibration on a Calabi-Yau manifold is essentially canonical when we fix the \k metric. This important property implies that one can construct the generalized special Lagrangian fibration over different part of the base $\partial \Delta$ separately and they will automatically match on the overlaps.\\ 

In this paper, we will construct generalized special Lagrangian fibration for $X_t$ over two types of regions in $\partial \Delta$. The type one region consists of the interior of top dimensional faces in $\partial \Delta$. Torus fibres over this region are of ``normal" size. This part of the construction is a direct analogue of Mclean's deformation result \cite{ML} and have also been done earlier by Goldstein in \cite{G} using a slightly different method. The type two region consists of small neighborhoods of vertices of $\partial \Delta$. Torus fibres over this region can have very ``small" size. This part of the construction is the main original contribution of this paper. Notice from \cite{lag3,toric} that these two types of regions are mirror dual to each other under the symplectic SYZ mirror symmetry proved in \cite{toric} for Calabi-Yau hypersurfaces in toric varieties.\\ 

All the resulting generalized special Lagrangian fibres are smooth tori. A (generalized) special Lagrangian fibration of smooth tori over an open set $U\subset \partial \Delta$ is said to {\bf represent the monodromy} if $\partial \Delta \setminus U$ is a fattening of $\tilde{\Gamma}$ that retracts to $\tilde{\Gamma}$. Existence of such monodromy representing (generalized) special Lagrangian (partial) fibration is a good indication that the actual global (generalized) special Lagrangian fibration might resemble the structure of the Lagrangian fibration constructed in \cite{lag1, lag2, lag3, toric} as we conjectured.\\

It is easy to see that the generalized special Lagrangian torus fibrations we construct over the union of the two types of regions already represent monodromy in the case of Fermat type Calabi-Yau hypersurfaces (including Fermat type quintics in $\mathbb{CP}^4$).\\

In the sequels of this paper, we will construct generalized special Lagrangian torus fibrations that represent monodromy for generic Calabi-Yau hypersurfaces near the large complex limit. We will also discuss some singular fibres.\\

This paper is organized as follows. In section 3 we discuss the implicit function theorems and the general deformation results of (generalized) special Lagrangians. Hamiltonian-gradient flow is reviewed in section 4, and is used to derive two more specific deformation results of (generalized) special Lagrangians that are needed for our construction. Sections 5 and 6 are devoted to the construction of generalized special Lagrangian fibrations over the two types of regions. In section 7, we apply our results to constructing monodromy representing generalized special Lagrangian torus fibrations for Fermat type Calabi-Yau hypersurfaces.\\

\se{Implicit function theorems and deformation of (generalized) special Lagrangian}
\subsection{Implicit function theorems}
In our discussion we will need the implicit function theorem (for example see \cite{GT}). We will carry out the standard proof carefully to derive a more quantitative version of the implicit function theorem that we need.\\

\begin{theorem}[Implicit function theorem]
\label{aa}
Assume that ${\cal B}_1$, ${\cal B}_2$ are Banach spaces. $F: {\cal B}_1\times \mathbb{R} \rightarrow {\cal B}_2$ is a map, whose Frechet derivative is continuous. $F(0,0)=0$ and $\frac{\partial F}{\partial h}(0,0): {\cal B}_1 \rightarrow {\cal B}_2$ is invertible. Then there exists a $C^1$-path $h(t)$ such that $F(h(t),t)=0$.
\end{theorem}

{\bf Proof:} Consider the operator $T$ defined as

\[
Th = h - \left(\frac{\partial F}{\partial h} (0,0)\right)^{-1} F(h,t).
\]

\[
\delta Th = \left(I - \left(\frac{\partial F}{\partial h} (0,0)\right)^{-1} \frac{\partial F}{\partial h}(h,t)\right) \delta h.
\]

$T$ would be a contraction mapping on the ball $U_{{\cal B}_1}(r) \subset {\cal B}_1$ if $\|[\frac{\partial F}{\partial h}(0,0)]^{-1}\| \leq C$ and 

\[
\left\|\frac{\partial F}{\partial h}(h,t) - \frac{\partial F}{\partial h}(0,0)\right\| \leq \frac{1}{2C}
\]

for $h \in U_{{\cal B}_1}(r)$. For $T$ to map $U_{{\cal B}_1}(r)$ to itself, we only need $\|F(0,t)\|_{{\cal B}_2} \leq \frac{r}{2C}$.\\

By choosing $r$ and $t_0$ small enough, all these conditions are easily satisfied for $(h,t) \in U_{{\cal B}_1}(r)\times [0,t_0]$. Then for any $t\in [0,t_0]$, by contraction mapping principle, there exists a unique $h(t) \in U_{{\cal B}_1}(r)$ satisfying $F(h(t),t)=0$. Clearly, $h(t)$ is a $C^1$-path with derivative

\[
\frac{dh}{dt} = - \left(\frac{\partial F}{\partial h}(h(t),t)\right)^{-1}\frac{\partial F}{\partial t}(h(t),t).
\]
\begin{flushright} \rule{2.1mm}{2.1mm} \end{flushright}
It is easy to see the above proof actually proved the following quantitative version of the implicit function theorem.

\begin{theorem}[Quantitative implicit function theorem]
\label{ab}
Assume that ${\cal B}_1$, ${\cal B}_2$ are Banach spaces. $F: {\cal B}_1\times \mathbb{R} \rightarrow {\cal B}_2$ is a map, whose Frechet derivative is continuous and $\frac{\partial F}{\partial h}(0,0): {\cal B}_1 \rightarrow {\cal B}_2$ is invertible. $\|[\frac{\partial F}{\partial h}(0,0)]^{-1}\| \leq C$. Assume further that we can find $r_0>r>0,t_0 >0$ such that for $(h,t) \in U_{{\cal B}_1}(r_0)\times [0,t_0]$,

\[
\left\|\frac{\partial F}{\partial h}(h,t) - \frac{\partial F}{\partial h}(0,0)\right\| \leq \frac{1}{2C},\ \ \|F(0,t)\|_{{\cal B}_2} \leq \frac{r}{2C}.
\]

Then there exists a $C^1$-path $h(t)$ in $U_{{\cal B}_1}(r)$ for $t\in [0,t_0]$ such that $F(h(t),t)=0$. Further more, $h(t)$ is the unique solution of $F(h(t),t)=0$ in $U_{{\cal B}_1}(r_0)$.
\end{theorem}
{\bf Proof:} Only the uniqueness part needs comment. The conditions of the theorem imply that $T: U_{{\cal B}_1}(r_0) \rightarrow U_{{\cal B}_1}(r_0)$ is a contraction mapping. Therefore $T$ has a unique fixed point in $U_{{\cal B}_1}(r_0)$, which is the unique solution of $F(h(t),t)=0$.
\begin{flushright} \rule{2.1mm}{2.1mm} \end{flushright}

\subsection{Deformation of generalized special Lagrangian}
Let $(X,\omega)$ be a symplectic manifold with a family of Calabi-Yau structures $(\Omega_t,J_t,g_t)$. Let $L_0$ be a generalized special Lagrangian submanifold in $(X,\omega,\Omega_0)$. Then locally, there is an identification $X\cong T^*L_0$ such that $\omega$ is identified with the canonical symplectic form on $T^*L_0$. (See \cite{W} for the general case. In our toric situation, such identification is actually explicit and does not need the general result from \cite{W}.) Any Lagrangian submanifold $L$ near $L_0$ can be identified with the graph of $dh$ for some smooth function $h$ on $L_0$. Under standard coordinate $(x,y)$ on $T^*L_0$, $\omega = dx\wedge dy$ and $L$ is locally the graph $(x, \frac{\partial h}{\partial x})$.\\

Locally, we can write

\[
\Omega_t = \eta_t\bigwedge_{k=1}^n (dx_k + u_{t,k}dy),
\]  
\[
\Omega_t|_L = \eta_t\left(x, \frac{\partial h}{\partial x}\right)\det\left(I + U_t\frac{\partial^2 h}{\partial x^2}\right)\bigwedge_{k=1}^n dx_k,
\] 

where $U_t = (u_{t,1},\cdots,u_{t,n}) = (u_{t,jk})$ and $u_{t,k}dy = \displaystyle\sum_{j=1}^n u_{t,jk}dy_j$.\\ 

\[
F(h,t) = {\rm Im}\left(\log \Omega_t|_L\right) = {\rm Im}\left(\log \eta_t\left(x, \frac{\partial h}{\partial x}\right)+ \log \det\left(I + U_t\frac{\partial^2 h}{\partial x^2}\right)\right)
\] 

defines a map $F: {\cal B}_1\times \mathbb{R} \rightarrow {\cal B}_2$, where ${\cal B}_1 = C^{2,\alpha}(L_0)$, ${\cal B}_2 = C^\alpha(L_0)$ are Banach spaces. It is also convenient to introduce ${\cal B}_3 = C^{1,\alpha}(\Gamma(T^*L_0))$ to measure Lagrangian submanifolds near $L_0$. ${\cal B}_1$ is naturally a closed subspace of ${\cal B}_3$. (More generally, a small tubular neighborhood of $L_0$ can be identified locally with the normal bundle $N_{L_0}$. A small deformation $L_1$ of $L_0$ can be viewed as a section of $N_{L_0}$. We will use $\|L_1 - L_0\|_{C^{1,\alpha}}$ to denote the $C^{1,\alpha}(\Gamma(N_{L_0}))$ norm of such section.) We intend to apply implicit function theorem to $F$ to construct the family of generalized special Lagrangians $L_t$ with respect to $(X,\omega,\Omega_t)$.\\

Straightforward computations give us

\[
\delta_h F(h,t) = \frac{\partial F}{\partial h}(h,t) \delta h = a_t^{ij}\frac{\partial^2 \delta h}{\partial x_i\partial x_j} + b_t^i\frac{\partial\delta h}{\partial x_i},
\]

where 

\[
a_t^{ij}(x) = {\rm Im}\left(\left(I + U_t\frac{\partial^2 h}{\partial x^2}\right)^{-1} U_t\right)\left(x,\frac{\partial h}{\partial x}\right),
\]
\[
b_t^i(x) = {\rm Im}\left(\frac{1}{\eta_t}\frac{\partial \eta_t}{\partial y_i} + {\rm Tr}\left(\left(I + U_t\frac{\partial^2 h}{\partial x^2}\right)^{-1}\frac{\partial U_t}{\partial y_i}\frac{\partial^2 h}{\partial x^2}\right)\right)\left(x,\frac{\partial h}{\partial x}\right).
\]

In particular,

\[
\delta_h F(0,0) = \frac{\partial F}{\partial h}(0,0) \delta h = a^{ij}\frac{\partial^2 \delta h}{\partial x_i\partial x_j} + b^i\frac{\partial\delta h}{\partial x_i},
\]

where 

\[
a^{ij}(x) = {\rm Im}\left((U_0)_{ij}\right)(x,0),\ \ b^i(x) = {\rm Im}\left(\frac{1}{\eta_0}\frac{\partial \eta_0}{\partial y_i}\right)(x,0).
\]

{\bf Example:} Assume that $(X,\omega,\Omega_0)$ is a toric variety with the \k potential $\rho$ (as function of $\{\log r_k^2\}$) and standard holomorphic volume form 

\[
\Omega_0 = \bigwedge_{k=1}^n \frac{-idz_k}{z_k} = \bigwedge_{k=1}^n (d\theta_k - i\frac{dr_k}{r_k}).
\]  

Then

\[
\omega = -i\partial\bar{\partial} \rho = -i\rho_{jk} \frac{dz_j}{z_j}\wedge \frac{d\bar{z}_k}{\bar{z}_k} = d\theta_k \wedge d\rho_k.
\]

Let $L_0$ be one of the real tori where all the $r_k$'s are constants. $L_0$ is obviously a generalized special Lagrangian with respect to the toric metric and $\Omega_0$. The coordinate $(\theta_k,\rho_k)$ naturally identifies $X$ with $T^*L_0$. Under this coordinate, it is easy to compute that

\[
\eta_0 =1,\ \ u_{0,jk} = \frac{-i}{2}\rho^{jk},\ \ {\rm where}\ (\rho^{jk}) = (\rho_{jk})^{-1}.
\]

Therefore

\begin{equation}
\label{ad}
\delta_h F(0,0) = \frac{-1}{2}\rho^{jk}\frac{\partial^2 \delta h}{\partial \theta_j\partial \theta_k}.
\end{equation}
\begin{flushright} \rule{2.1mm}{2.1mm} \end{flushright}
A smooth function is called (smoothly) bounded if the function and all its multi-derivatives are bounded. The metric matrix $(g_{ij})$ under certain coordinate is called bounded if each entry $g_{ij}$ is bounded as a smooth function and $\det(g_{ij})$ is bounded from below by a positive constant. Under such condition, the curvature will be smoothly bounded and the inverse matrix $(g^{ij})$ is also bounded in the same sense. The Riemannian manifold is of bounded geometry if in addition, the injective radius is bounded from below by a positive constant. Such coordinates that make $(g_{ij})$ bounded is called proper.\\ 

A smooth family of smooth embeddings $\{\phi_t: Y_1 \rightarrow (Y_2,g)\}_{t\in B}$ is called bounded if $y_2 = \phi_t(y_1) = \phi(y_1,t)$ is smoothly bounded for suitable coordinate $(y_1,t)$ and $y_2$ such that $y_1$ for fixed $t$ is proper coordinate of $(Y_1,\phi_t^*g)$ and $y_2$ is proper coordinate of $(Y_2,g)$. (Smoothly) boundedness of a smooth family of smooth tensors can also be defined similarly.\\ 

\begin{de}
$(L \subset X ,g)$ is said to be of bounded geometry if both $(X,g)$ and $(L,g|_L)$ are of bounded geometry and the embedding $L \rightarrow (X,g)$ is bounded.\\
\end{de}

\begin{theorem}
\label{ac}
Let $(X,\omega)$ be a symplectic manifold with a family of Calabi-Yau structure $(\Omega_t,J_t,g_t)$. Let $L_0$ be a generalized special Lagrangian submanifold in $(X,\omega,\Omega_0)$. Assume $(L_0 \subset X, g_0)$ is of bounded geometry and $\tilde{\nu} = \max_{t\in [0,t_0]}\|\Omega_t - \Omega_0\|_{C^2}$ is small enough. Then there exist $C_1,C_2>0$ and a smooth family of generalized special Lagrangian submanifolds $L_t \subset (X,\omega,\Omega_t)$ that are Hamiltonian equivalent to $L_0$ for $t\in [0,t_0]$, and satisfy $L_t \in U_{{\cal B}_3}(C_1\tilde{\nu})$. Further more, $L_t$ is the unique generalized special Lagrangian submanifold of $(X,\omega,\Omega_t)$ in $U_{{\cal B}_3}(C_2)$.
\end{theorem}
{\bf Proof:} $\frac{\partial F}{\partial h}(h,t)$ is a second order linear differential operator with coefficients being smooth functions of $(t,x,h,\frac{\partial h}{\partial x},\frac{\partial^2 h}{\partial x^2})$. Therefore $\frac{\partial F}{\partial h}(h,t)$ is continuous on $(h,t)$. Similarly, $\frac{\partial F}{\partial t}(h,t)$ is continuous on $(h,t)$.\\ 

In general, as is well known, $J_t$ is determined by $\Omega_t$ according to $J_t(dx_k + {\rm Re}(u_{t,k}dy)) = {\rm Im}(u_{t,k}dy)$, which easily implies that ${\rm Im}(U_t)(x,0) = (g_t|_L)^{-1}$. In particular

\[
\frac{\partial F}{\partial h}(0,0) \delta h = (g_0|_L)^{ij}\frac{\partial^2 \delta h}{\partial x_i\partial x_j} + b^i\frac{\partial\delta h}{\partial x_i}.
\]

(According to (\ref{ad}), the operator is much simpler in our toric case.) Since $(L_0 \subset X, g_0)$ is of bounded geometry, by standard Schauder estimate for the linear elliptic operator $\frac{\partial F}{\partial h}(0,0)$, there exists a constant $C$ such that

\[
\|\delta h\|_{{\cal B}_1} \leq C \left\|\left(\frac{\partial F}{\partial h}(0,0)\right)\delta h\right\|_{{\cal B}_2}.
\]

Since $[0,t_0]$ is compact and $\frac{\partial F}{\partial h}(h,t)$ is continuous on $(h,t)$, there exists a constant $C_2$ such that when $(h,t) \in U_{{\cal B}_1}(C_2)\times [0,t_0]$,

\[
\left\|\frac{\partial F}{\partial h}(h,t) - \frac{\partial F}{\partial h}(0,t)\right\| \leq \frac{1}{4C}.
\]

(More quantitatively, one needs the $C^2$-norm of $\Omega_t$ to be uniformly bounded, which is implied by the conditions of the theorem.)\\

The condition $\tilde{\nu} = \max_{t\in [0,t_0]}\|\Omega_t - \Omega_0\|_{C^2}$ being small implies both 

\[
\left\|\frac{\partial F}{\partial h}(0,t) - \frac{\partial F}{\partial h}(0,0)\right\| \leq \frac{1}{4C}
\]

and

\[
\|F(0,t)\|_{{\cal B}_2} \leq \frac{C_1\tilde{\nu}}{2C}.
\]

Now our quantitative implicit function theorem will imply our result.
\begin{flushright} \rule{2.1mm}{2.1mm} \end{flushright}
{\bf Remark:} Since $\|\Omega_t - \Omega_0\|_{C^2} = O(t)$ for $t$ small, when $t_0$ is allowed to be small, no additional condition on the smooth family $\{\Omega_t\}_{t\in [0,t_0]}$ is needed. Theorem in this version is an analog of Mclean's theorem.\\

\se{Hamiltonian-gradient flow}
Consider a smooth family of \k manifold $(X_t,\omega_t)$ with constant \k class, total space ${\cal X}$ and fibration map $t: {\cal X} \rightarrow D \subset \mathbb{C}$. It is very easy to construct a \k form $\omega$ on the total space ${\cal X}$ so that $\omega_t = \omega|_{X_t}$. $\omega$ is determined by the family up to the pullback of \k form on $D$. The flow of $V=\frac{\nabla f}{|\nabla f|^2}$, where $f = {\rm Re}(t)$, determines symplectomorphisms $\phi_t: X_0 \rightarrow X_t$. We call this flow Hamiltonian-gradient flow. $V$ will be called the (normalized) Hamiltonian-gradient vector field for the parameter $t$. This flow was discussed intensively in \cite{lag2} and has been used in my construction of Lagrangian fibrations \cite{lag1, lag2, lag3, toric}.\\

We will start with a family $(X_t,\omega_t,\Omega_t)$ of Calabi-Yau manifolds, we have the total space $({\cal X},\omega,\Omega)$, such that $\omega_t = \omega|_{X_t}$ and $\Omega_t = (\Omega\otimes (dt)^{-1})|_{X_t}$. Since $\phi_t$ are symplectomorphisms, we have $\omega_0 = \phi_t^*\omega_t$. Using the Hamiltonian-gradient flow, the family $(X_t,\omega_t,\Omega_t)$ of Calabi-Yau manifolds can be reduced to the equivalent family of Calabi-Yau structures $\{\phi_t^*\Omega_t\}$ on the fixed symplectic manifold $(X_0,\omega_0)$. Let $V_t$ denote the Hamiltonian-gradient vector field for the parameter $t$. We have\\ 

\begin{lm}
\[
\imath(V_t)\Omega|_{X_t} = \Omega_t.
\]
\end{lm}

{\bf Proof:} By definition, $\langle V_t,dt\rangle =1$. Hence $\imath(V_t)\Omega$ is a representative of the relative holomorphic volume form, namely

\[
\imath(V_t)\Omega|_{X_t} = \Omega_t.
\]
\begin{flushright} \rule{2.1mm}{2.1mm} \end{flushright}
\begin{theorem}
\label{da}
Assume that $\{X_t\}_{t\in [0,t_0]}$ is a bounded smooth family, and the bound is uniform with respect to $t_0$. $(L_0\subset X_0,g_0)$ is of bounded geometry. Then when $t_0$ is small, there exist a smooth family of generalized special Lagrangian submanifolds $L_t \subset (X_t,\omega_t,\Omega_t)$ that are Hamiltonian equivalent to $L_0$ for $t\in [0,t_0]$, and satisfy $\|L_t - L_0\|_{C^{1,\alpha}} \leq C|t|$.
\end{theorem}
{\bf Proof:} Under the condition of the theorem, there exist local bounded smooth map, $(z,t) \in \mathbb{C}^n\times [0,t_0] \rightarrow {\mathcal X}$, such that $(z,t)$ for fixed $t$ is a local proper coordinate of $X_t$. The projection of $\frac{\partial}{\partial t}$ to the normal direction of $X_t$ is exactly the Hamiltonian-gradient vector field $V_t$. $V_t$ is clearly smooth and uniformly bounded under our condition. Consequently, $\phi_t$ is a uniformly bounded smooth flow, $\Omega_t = \imath(V_t)\Omega|_{X_t}$ and $\phi_t^*\Omega_t$ are smooth uniformly bounded family. Therefore $\|\phi_t^*\Omega_t - \Omega_0\|_{C^2} = O(t)$. For $t_0$ small, apply theorem \ref{ac}, there exists a smooth Hamiltonian equivalent family of generalized special Lagrangian submanifolds $L'_t \subset (X_0,\omega_0,\phi_t^*\Omega_t)$ for $t\in [0,t_0]$. Let $L_t = \phi_t(L'_t)$. Clearly, $\{L_t\}_{t\in[0,t_0]}$ is the desired smooth family.
\begin{flushright} \rule{2.1mm}{2.1mm} \end{flushright}
We will also need to consider a family version (parametrized by $s$) of the above situation, where we have 2-parameter family $(X_{t,s},\omega_{t,s},\Omega_t)$ of Calabi-Yau hypersurfaces in $({\cal X},\omega,\Omega)$, such that $\omega_{t,s} = \omega|_{X_{t,s}}$ and $\Omega_t = (\Omega\otimes (dt)^{-1})|_{X_{t,s}}$. (Here $\Omega$, $\Omega_t$ are also depending on $s$. We are omitting the subscript for the simplicity of notations.) We will use $V_t$ and $V$ to denote Hamiltonian-gradient vector fields for the parameters $t$ and $s$ respectively. We still have $\imath(V_t)\Omega|_{X_{t,s}} = \Omega_t$. Let $\phi_s$ denote the flow for $V$. For fixed $t$, the family $(X_{t,s},\omega_{t,s},\Omega_t)$ of Calabi-Yau manifolds can be reduced to the equivalent family of Calabi-Yau structures $\{\phi_s^*\Omega_t\}$ on the fixed symplectic manifold $(X_{t,0},\omega_{t,0})$.

\begin{theorem}
\label{db}
Assume that $X_{\tilde{t},s}$ forms a bounded smooth family for $\tilde{t}$ near $t$ and $s\in [0,s_0]$, and the bound is uniform with respect to $s_0$. $(L_{t,0}\subset X_{t,0},g_0)$ is of bounded geometry. Then when $s_0$ is small, there exist a smooth family of generalized special Lagrangian submanifolds $L_{t,s} \subset (X_{t,s},\omega_s,\Omega_s)$ that are Hamiltonian equivalent to $L_{t,0}$ for $s\in [0,s_0]$, and satisfy $\|L_{t,s} - L_{t,0}\|_{C^{1,\alpha}} \leq C|s|$.
\end{theorem}
{\bf Proof:} Under the condition of the theorem, there exists a local bounded smooth map, $(z,\tilde{t},s) \in \mathbb{C}^n\times (t-\epsilon,t+\epsilon)\times [0,s_0] \rightarrow {\mathcal X}$, such that $(z,\tilde{t},s)$ for fixed $(\tilde{t},s)$ is a local proper coordinate of $X_{\tilde{t},s}$. The projection of $\frac{\partial}{\partial \tilde{t}}$ and $\frac{\partial}{\partial s}$ to the normal direction of $X_{t,s}$ is exactly the Hamiltonian-gradient vector fields $V_t$ and $V$. $V_t$ and $V$ are clearly smooth and uniformly bounded under our condition. Consequently, $\phi_s$ is a uniformly bounded smooth flow, $\Omega_s = \imath(V_t)\Omega|_{X_{t,s}}$ and $\phi_s^*\Omega_s$ are smooth uniformly bounded family. Therefore $\|\phi_s^*\Omega_s - \Omega_0\|_{C^2} = O(t)$. For $s_0$ small, applying theorem \ref{ac}, there exists a smooth Hamiltonian equivalent family of generalized special Lagrangian submanifolds $L'_{t,s} \subset (X_{t,0},\omega_0,\phi_s^*\Omega_s)$ for $s\in [0,s_0]$. Let $L_{t,s} = \phi_s(L'_{t,s})$. Clearly, $\{L_{t,s}\}_{s\in[0,s_0]}$ is the desired smooth family.
\begin{flushright} \rule{2.1mm}{2.1mm} \end{flushright}

\se{Fibration over top dimensional faces in $\partial \Delta$}
Take $X_0 \subset P_\Delta$ to be the large complex limit, with the natural generalized special Lagrangian fibration $\hat{\pi}_0: X_0 \rightarrow \partial \Delta$. Non-degenerate fibres are in the smooth part of $X_0$, which is a union of top dimensional complex tori fibred over top dimensional faces in $\partial \Delta$. The following results are direct analogues of Mclean's deformation result \cite{ML} and theorem \ref{eb} in a slightly different form has also been done earlier by Goldstein in \cite{G} using a slightly different method.\\ 

\begin{theorem}
\label{eb}
For any non-degenerate generalized special Lagrangian fibre $L_0$ in $X_0$, there exists a smooth family $\{L_t\}_{t\in[0,t_0]}$, where $L_t$ is a generalized special Lagrangian torus in $X_t$ that is Hamiltonian equivalent to $\phi_t(L_0)$ and satisfy $\|L_t - L_0\|_{C^{1,\alpha}} \leq C|t|$. ($C$ and $t_0$ depend on the distance of $L_0$ to ${\rm Sing}(X_0)$). For each fixed small $t$, when the fibre $L_0 \subset X_0$ varies, $L_t \subset X_t$ will form a generalized special Lagrangian fibration over an open set $U_t^{\rm top} \subset \partial \Delta$. The compliment of $U_t^{\rm top}$ in $\partial \Delta$ is a thin fattening (depending how small the $t$ is) of the $(n-1)$-skeleton of $\partial \Delta$.
\end{theorem}

{\bf Proof:} Without loss of generality, we will concentrate on one of the smooth components of $X_0$. One may choose suitable toric coordinate $(z_0,\cdots,z_n)$, such that locally $X_0 = \{z_0=0\}$ and $X_t = \{z_0 = tp(z)\}$, where $p(z)$ is holomorphic on this smooth component of $X_0$. A generalized special Lagrangian fibre in this smooth component of $X_0$ can be expressed as

\begin{equation}
\label{ea}
L_0 = \{z_0=0, |z_k|=r_k ({\rm constant}),\ {\rm for}\ 1\leq k \leq n\}.
\end{equation}

Since $L_0$ is non-degenerate, $p(z)$ is bounded near $L_0$. Hence $\{X_t\}_{t\in [0,t_0]}$ is a smooth family that is bounded near $L_0$ (uniformly with respect to $t_0$ as long as $t_0$ is not too big). Applying theorem \ref{da}, we get the desired smooth family $\{L_t\}_{t\in[0,t_0]}$.\\

For fixed small $t$, to show $L_t$ form a generalized special Lagrangian fibration when the fibre $L_0$ varies, it is sufficient to show that non-trivial deformation 1-forms of $L_t$ have no zeroes on $L_t$. This is true because $\|L_t - L_0\|_{C^{1,\alpha}} \leq C|t|$ is small, which implies that non-trivial deformation 1-forms of $L_t$ are close to non-trivial harmonic 1-forms of $L_0$, which have no zeroes on $L_0$. (Notice from (\ref{ad}) that for $L_0$, deformation 1-forms are just harmonic 1-forms.)
\begin{flushright} \rule{2.1mm}{2.1mm} \end{flushright}
We will call $L_0$ in (\ref{ea}) of bounded geometry up to scale $\nu$ if $\log \hat{r}_k$ are bounded for $1\leq k \leq n$, where $\hat{r}_k = r_k/\nu$. Argument in theorem \ref{eb} will also apply to such torus as long as $\hat{t} = t/\nu^{n+1}$ is small. More precisely, under the rescaled metric $\hat{g} = \nu^{-2}g$ and coordinate $\hat{z} = z/\nu$ and $\hat{t}$, we are exactly in the situation of the proof of theorem \ref{eb}, where we can apply theorem \ref{da}. We will use $\hat{C}^{1,\alpha}$ to denote $C^{1,\alpha}$-norm with respect the scalled metric $\hat{g}$.\\

\begin{theorem}
\label{ec}
$U_t^{\rm top}$ in theorem \ref{eb} will contain $\hat{\pi}_0(L_0)$, where $L_0 \subset X_0$ is of bounded geometry up to scale $\nu$ as long as $\hat{t} = t/\nu^{n+1}$ is small. Further more $\|L_t - L_0\|_{\hat{C}^{1,\alpha}} = O(\hat{t})$ for such torus.
\end{theorem}
\begin{flushright} \rule{2.1mm}{2.1mm} \end{flushright}

\se{Fibration near vertices of $\partial \Delta$}
The fibre of $\hat{\pi}_0: X_0 \rightarrow \partial \Delta$ over a vertex of $\Delta$ is a fixed point in $P_\Delta$ of the toric action. Assume $O\in P_\Delta$ is one of such fixed points and $P_\Delta$ is smooth around $O$. Then locally around $O$, we have local toric coordinate $z=(z_0,\cdots, z_n)$ such that

\[
X_t = \{\tilde{p}_t(z)=0\},\ \ {\rm where}\ \tilde{p}_t(z) = \prod_{k=0}^n z_k + tp(z).
\]

A natural idea to construct generalized special Lagrangian for $X_t$ (which is similar to the idea we used on the construction of Lagrangian fibration in \cite{lag1, lag2, lag3, toric}) is to modify the singular Hamiltonian-gradient flow to deform the generalized special Lagrangian for $X_0$ to the generalized special Lagrangian for $X_t$. This idea was successfully carried out in the previous section for smooth fibres in $X_0$. For singular fibres of $X_0$, this idea is somewhat difficult to carry out, partly due to the singular nature of the Hamiltonian-gradient vector field.\\

Instead, we will use an alternative idea based on a local model so that we will not need to deal with singular vector fields. We will start with the local model.\\

{\bf Local model:} In $\mathbb{C}^{n+1}$, consider the family of hypersurfaces $\{X_{t,0}\}$ defined as

\[
X_{t,0} = \{z_0\cdots z_n =t\}.
\]

Clearly $X_{t,0}$ are all toric varieties. For a toric metric on $\mathbb{C}^{n+1}$ with \k form $\omega$, $\omega_{t,0} = \omega|_{X_{t,0}}$ is a toric \k form on $X_{t,0}$. The natural generalized special Lagrangian torus fibration $\hat{\pi}_{t,0}: X_{t,0} \rightarrow \partial \Delta$ with respect to any toric \k form is defined as $\hat{\pi}_{t,0}(z) = \{|z_k|^2-\min(\{|z_i|^2\}_{i=0}^n)\}_{k=0}^n$, where $\Delta$ is the first quadrant in $\mathbb{R}^{n+1}$.
\begin{flushright} \rule{2.1mm}{2.1mm} \end{flushright}
Without loss of generality, we may assume $p(0)=1$ and denote $\check{p}(z) = p(z) - 1$. Let

\[
X_{t,s} = \{\tilde{p}_{t,s}(z)=0\},\ \ {\rm where}\ \tilde{p}_{t,s}(z) = \prod_{k=0}^n z_k + t(1 + s\check{p}(z)).
\]

We will use the family $\{X_{t,s}\}_{s\in [0,1]}$, which connect the local model $X_{t,0}$ and $X_t = X_{t,1}$. Our idea here is to deform the generalized special Lagrangian torus fibration of the local model $X_{t,0}$ to the generalized special Lagrangian torus fibration of $X_t$.\\

For a constant $\nu>0$, let

\[
L_{t,0} = \{z\in X_{t,0}||z_k|=r_k ({\rm constant}),\ {\rm for}\ 0\leq k \leq n\},
\]

where $r_k$ for $1\leq k \leq n$ are of order $\nu$ and $r_0$ is of order $\frac{t}{\nu^n}$. More precisely, assume that $\log \frac{r_k}{\nu}$ for $1\leq k \leq n$ and $\log \frac{\nu^n r_0}{t}$ are bounded. Similar to section 5, we will use rescaled metric $\hat{g} = \nu^{-2}g$ and coordinate $\hat{z} = z/\nu$ and $\hat{t} = t/\nu^{n+1}$. Notice that when $\hat{t}$ is bounded, $L_{t,0}$ is of bounded geometry up to scale $\nu$. We have\\

\begin{theorem}
\label{fa}
When $\hat{t}$ and $\nu$ are bounded and $t$ is small enough, for any such generalized special Lagrangian fibre $L_{t,0}$ in $X_{t,0}$, there exist a smooth family $\{L_{t,s}\}_{s\in[0,1]}$, where $L_{t,s}$ is a generalized special Lagrangian torus in $X_{t,s}$ that is Hamiltonian equivalent to $\phi_s(L_{t,0})$ and $\|L_{t,s} - L_{t,0}\|_{\hat{C}^{1,\alpha}} = O(\hat{t}\nu)$. For each fixed small $t$, when the fibre $L_{t,0}$ varies, $L_{t,1}$ will form a generalized special Lagrangian fibration over an open neighborhood $U_t^{\rm ver}$ of the vertex of $\partial \Delta$.
\end{theorem}
{\bf Proof:} Scale the coordinate and metric by $z = \nu \hat{z}$, $g = \nu^2 \hat{g}$. Then 

\begin{equation}
\label{fb}
X_{t,s} = \left\{\hat{z}_0 = - \hat{t}(1+s\check{p}(\nu\hat{z}))\left(\prod_{k=1}^n \hat{z}_k\right)\right\}.
\end{equation}

Near $L_{t,0}$, $\hat{z}_k = O(1)$ for $1\leq k \leq n$ and $\hat{z}_0 = O(\hat{t})$, which can be small. From the explicit expression in (\ref{fb}), it is easy to verify that $X_{t,s}$ (near $L_{t,0}$) with respect to parameters $\hat{t}$ near the fixed value and $\hat{s} = \hat{t}\nu s\in [0,\hat{t}\nu]$ forms a smooth family that is uniformly bounded. According to our requirement on $L_{t,0}$, clearly $L_{t,0}$ is of bounded geometry in $(X_{t,0},\hat{g})$. Applying theorem \ref{db}, we get the desired smooth family $\{L_{t,s}\}_{s\in[0,1]}$.
\begin{flushright} \rule{2.1mm}{2.1mm} \end{flushright}

\se{Monodromy representing generalized special Lagrangian torus fibration for Fermat type Calabi-Yau hypersurfaces}
The results in the previous two sections can be summarized into the following.\\  

\begin{theorem}
\label{ga}
When $t$ is small enough, there exists an open subset $\tilde{X}_t \subset X_t$ and a smooth special Lagrangian torus fibration $\hat{\pi}_t: \tilde{X}_t \rightarrow U_t$, where $U_t = U_t^{\rm top} \cup U_t^{\rm ver} \subset \partial \Delta$ (as indicated in figure 1 for $n=2$ near a vertex of $\Delta$).
\end{theorem}
{\bf Proof:} The only thing remaining to be shown is that the two fibrations coincide in the overlaps. Start with

\[
L_{t,0} = \{z\in X_{t,0}||z_k|=r_k ({\rm constant}),\ {\rm for}\ 0\leq k \leq n\}
\]
\[
L_0 = \{z_0=0, |z_k|=r_k ({\rm constant}),\ {\rm for}\ 1\leq k \leq n\}
\]

with the same constants $\{r_k\}_{k=1}^n$ such that $L_0$ is of bounded geometry up to scale $\nu$. Under the rescaling by $\nu$, it is easy to check that $\|L_{t,0} - L_0\|_{\hat{C}^{1,\alpha}} = O(\hat{t})$. By theorems \ref{ec} and \ref{fa}, we have

\[
\|L_t - L_0\|_{\hat{C}^{1,\alpha}} = O(\hat{t}),\ \ \|L_{t,1} - L_{t,0}\|_{\hat{C}^{1,\alpha}} = O(\hat{t}).
\]

Consequently $\|L_{t,1} - L_t\|_{\hat{C}^{1,\alpha}} = O(\hat{t})$. This implies that the Lagrangian torus $L_{t,1} \subset X_t$ is of order $O(\hat{t})$ in the deformation space of Lagrangian torus (modulo Hamiltonian equivalence) near $L_t \subset X_t$ parametrized by $H^1(L_t)$ under the rescaled metric $\hat{g}$. Recall that $L_t$ form a generalized special Lagrangian fibration with bounded geometry with respect to $\hat{g}$ when constants $\{r_k\}_{k=1}^n$ vary. We can find constants $\{r'_k\}_{k=1}^n$ such that $\log r'_k - \log r_k = O(\hat{t})$ and the corresponding $L'_t$ is Hamiltonian equivalent to $L_{t,1}$. It is easy to verify from their exlicit definition that $\|L_{t,0} - L'_0\|_{\hat{C}^{1,\alpha}} \leq \|L_{t,0} - L_0\|_{\hat{C}^{1,\alpha}} + \|L_0 - L'_0\|_{\hat{C}^{1,\alpha}} = O(\hat{t})$. Consequently, $\|L_{t,1} - L'_t\|_{\hat{C}^{1,\alpha}} \leq \|L_{t,1} - L_{t,0}\|_{\hat{C}^{1,\alpha}} + \|L_{t,0} - L'_0\|_{\hat{C}^{1,\alpha}} + \|L'_0 - L'_t\|_{\hat{C}^{1,\alpha}} = O(\hat{t})$. By the uniqueness part of theorem \ref{ac}, when $\hat{t}$ is small, $L_{t,1}$ coincide with $L'_t$.
\begin{flushright} \rule{2.1mm}{2.1mm} \end{flushright}
\begin{center}
\begin{picture}(200,170)(50,20)
\thinlines
\multiput(82,149)(130,0){2}{\vector(2,1){38}}
\multiput(82,149)(130,0){2}{\vector(-2,1){38}}
\multiput(82,149)(130,0){2}{\vector(0,-1){38}}
\multiput(119,160)(130,0){2}{\footnotesize{$|z_1|$}}
\multiput(33,160)(130,0){2}{\footnotesize{$|z_2|$}}
\multiput(84,111)(130,0){2}{\footnotesize{$|z_0|$}}
\thicklines
\multiput(82,149)(130,0){2}{\circle{15}}
\put(82,149){\circle{38}}
\put(90,149){\line(1,0){12}}
\put(74,149){\line(-1,0){12}}
\put(86,155){\line(3,5){7}}
\put(78,143){\line(-3,-5){7}}
\put(78,155){\line(-3,5){7}}
\put(86,143){\line(3,-5){7}}
\multiput(77,168)(3,0){4}{\circle*{1}}
\multiput(74,165)(3,0){6}{\circle*{1}}
\multiput(76,162)(3,0){5}{\circle*{1}}
\multiput(78,159)(3,0){4}{\circle*{1}}
\multiput(79,156)(3,0){3}{\circle*{1}}
\multiput(78,153)(3,0){4}{\circle*{1}}
\multiput(76,150)(3,0){5}{\circle*{1}}
\multiput(64,147)(3,0){13}{\circle*{1}}
\multiput(64,144)(3,0){13}{\circle*{1}}
\multiput(66,141)(3,0){4}{\circle*{1}}
\multiput(69,138)(3,0){2}{\circle*{1}}
\multiput(71,135)(3,0){1}{\circle*{1}}
\multiput(90,141)(3,0){4}{\circle*{1}}
\multiput(93,138)(3,0){2}{\circle*{1}}
\multiput(94,135)(3,0){1}{\circle*{1}}

\put(220,149){\line(1,0){24}}
\put(204,149){\line(-1,0){24}}
\put(216,155){\line(3,5){14}}
\put(208,143){\line(-3,-5){14}}
\put(208,155){\line(-3,5){14}}
\put(216,143){\line(3,-5){14}}

\multiput(197,179)(3,0){12}{\circle*{1}}
\multiput(199,176)(3,0){11}{\circle*{1}}
\multiput(200,173)(3,0){9}{\circle*{1}}
\multiput(202,170)(3,0){8}{\circle*{1}}
\multiput(204,167)(3,0){7}{\circle*{1}}
\multiput(206,164)(3,0){5}{\circle*{1}}
\multiput(208,161)(3,0){4}{\circle*{1}}
\multiput(209,158)(3,0){3}{\circle*{1}}

\multiput(182,147)(3,0){8}{\circle*{1}}
\multiput(183,144)(3,0){8}{\circle*{1}}
\multiput(185,141)(3,0){8}{\circle*{1}}
\multiput(186,138)(3,0){7}{\circle*{1}}
\multiput(188,135)(3,0){6}{\circle*{1}}
\multiput(189,132)(3,0){5}{\circle*{1}}
\multiput(191,129)(3,0){3}{\circle*{1}}
\multiput(192,126)(3,0){2}{\circle*{1}}
\multiput(193,123)(3,0){1}{\circle*{1}}

\multiput(222,147)(3,0){8}{\circle*{1}}
\multiput(220,144)(3,0){8}{\circle*{1}}
\multiput(219,141)(3,0){8}{\circle*{1}}
\multiput(221,138)(3,0){7}{\circle*{1}}
\multiput(223,135)(3,0){6}{\circle*{1}}
\multiput(226,132)(3,0){4}{\circle*{1}}
\multiput(228,129)(3,0){3}{\circle*{1}}
\multiput(229,126)(3,0){2}{\circle*{1}}
\multiput(231,123)(3,0){1}{\circle*{1}}

\thicklines
\put(140,120){\line(2,1){12}}
\put(140,120){\line(-2,1){12}}
\put(140,120){\vector(0,-1){12}}
\put(112,125){$U_t^{\rm ver}$}
\put(155,125){$U_t^{\rm top}$}
\put(137,98){$U_t$}

\thinlines
\put(142,59){\vector(2,1){38}}
\put(142,59){\vector(-2,1){38}}
\put(142,59){\vector(0,-1){38}}
\put(179,70){\footnotesize{$|z_1|$}}
\put(93,70){\footnotesize{$|z_2|$}}
\put(144,21){\footnotesize{$|z_0|$}}
\thicklines
\put(142,59){\circle{15}}
\multiput(139,66)(3,0){3}{\circle*{1}}
\multiput(138,63)(3,0){4}{\circle*{1}}
\multiput(136,60)(3,0){5}{\circle*{1}}

\put(150,59){\line(1,0){24}}
\put(134,59){\line(-1,0){24}}
\put(146,65){\line(3,5){14}}
\put(138,53){\line(-3,-5){14}}
\put(138,65){\line(-3,5){14}}
\put(146,53){\line(3,-5){14}}

\multiput(127,89)(3,0){12}{\circle*{1}}
\multiput(129,86)(3,0){11}{\circle*{1}}
\multiput(130,83)(3,0){9}{\circle*{1}}
\multiput(132,80)(3,0){8}{\circle*{1}}
\multiput(134,77)(3,0){7}{\circle*{1}}
\multiput(136,74)(3,0){5}{\circle*{1}}
\multiput(138,71)(3,0){4}{\circle*{1}}
\multiput(139,68)(3,0){3}{\circle*{1}}

\multiput(112,57)(3,0){13}{\circle*{1}}
\multiput(113,54)(3,0){12}{\circle*{1}}
\multiput(115,51)(3,0){8}{\circle*{1}}
\multiput(116,48)(3,0){7}{\circle*{1}}
\multiput(118,45)(3,0){6}{\circle*{1}}
\multiput(119,42)(3,0){5}{\circle*{1}}
\multiput(121,39)(3,0){3}{\circle*{1}}
\multiput(122,36)(3,0){2}{\circle*{1}}
\multiput(123,33)(3,0){1}{\circle*{1}}

\multiput(152,57)(3,0){8}{\circle*{1}}
\multiput(150,54)(3,0){8}{\circle*{1}}
\multiput(149,51)(3,0){8}{\circle*{1}}
\multiput(151,48)(3,0){7}{\circle*{1}}
\multiput(153,45)(3,0){6}{\circle*{1}}
\multiput(156,42)(3,0){4}{\circle*{1}}
\multiput(158,39)(3,0){3}{\circle*{1}}
\multiput(159,36)(3,0){2}{\circle*{1}}
\multiput(161,33)(3,0){1}{\circle*{1}}
\end{picture}
\end{center}
\begin{center}
\stepcounter{figure}
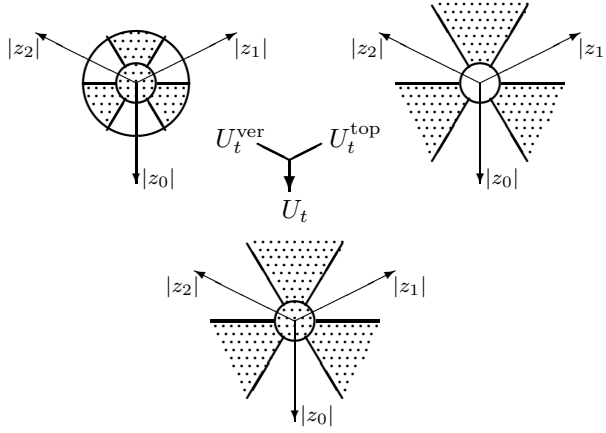
Figure \thefigure: $U_t = U_t^{\rm top} \cup U_t^{\rm ver} \subset \partial \Delta$  (n=2)
\end{center}
The dotted areas in figure 1 are $U_t^{\rm ver}$, $U_t^{\rm top}$ and $U_t = U_t^{\rm top} \cup U_t^{\rm ver}$ near a vertex of $\partial \Delta$ (for $n=2$). The radii of the small circles are of order $O(|t|^{\frac{1}{n+1}})$ and the radius of the big circle is of order $O(1)$, which is independent of $t$.\\

{\bf Remark:} By ``smooth special Lagrangian torus fibration", we mean that the fibres are smooth and varying smoothly. Then $U_t$ can be equipped with the corresponding smooth structure to make the map $\hat{\pi}_t$ smooth in the usual sense.\\

When $\dim X_t =1$, the singular locus $\tilde{\Gamma} = \emptyset$. Theorem \ref{ga} implies the uniformization theorem for such $X_t$, which is of course classical.\\

\begin{co}
When $\dim X_t =1$ and $t$ is small enough, there exists a smooth special Lagrangian circle fibration $\hat{\pi}_t: X_t \rightarrow S^1$, consequently, $X_t$ is biholomorphic to $\mathbb{C}/\Pi$, where $\Pi$ is a lattice in $\mathbb{C}$.\\
\end{co}
{\bf Proof:} Assume the phase of the holomorphic volume form $\Omega_t$ is $\theta$ along the special Lagrangian circle fibres. One can find canonical (up to addition of constant) local holomorphic coordinate $z=x+iy$ such that $\Omega_t = e^{i\theta}dz$, and $dy$ resticted to the special Lagrangian circle fibres equals to zero. Consequently, $y$ restricts to constant along the special Lagrangian circles and forms natual coordinate of the base $\partial \Delta \cong S^1$. When extending the canonical holomorphic coordinate $z$ around the fibre circle and a horizontal circle, one would pick up constants $\pi_1 \in \mathbb{R}$ and $\pi_2 \in \mathbb{C}$ that generate lattice $\Pi\subset\mathbb{C}$. Then it is easy to see that $X_t$ is biholomorphic to $\mathbb{C}/\Pi$ and the fibration map is $\hat{\pi}_t(z) = y$.
\begin{flushright} \rule{2.1mm}{2.1mm} \end{flushright}
A more non-trival application of theorem \ref{ga} is to the case of Fermat type Calabi-Yau in $\mathbb{P}^n$ (including Fermat type quintics in $\mathbb{P}^4$). For Fermat type Calabi-Yau, the singular locus $\tilde{\Gamma}$ is rather simple. (See \cite{lag1,lag2} for detail.) It is fairly easy to see that $\partial \Delta \setminus U_t$ is a fattening of $\tilde{\Gamma}$ in this case.\\

\begin{co}
For Fermat type Calabi-Yau family $\{X_t\}$ with the restriction of any toric \k metric from $\mathbb{P}^n$ (including the Fubini-Study metric), when $t$ is small enough, there exists a smooth monodromy representing generalized special Lagrangian torus (partial) fibration for $X_t$.
\end{co}
\begin{flushright} \rule{2.1mm}{2.1mm} \end{flushright}
{\bf Acknowledgement:} I would like to thank Prof. Yong-Geun Oh for helpful discussions and Qin Jing for pushing me to draw figure 1.\\

\ifx\undefined\bysame
\newcommand{\bysame}{\leavevmode\hbox to3em{\hrulefill}\,}
\fi

\end{document}